\documentclass[a4paper,11pt,reqno]{amsart}
\input{glyphtounicode}
\makeatletter
\@ifundefined{subjclassname@2020}{%
  \@namedef{subjclassname@2020}{\textup{2020} Mathematics Subject Classification}}{}
\makeatother
\usepackage{amssymb}
\usepackage{ifthen}
\usepackage{graphicx}
\usepackage{float}
\usepackage{caption}
\usepackage{subcaption}
\usepackage{cite}
\usepackage{amsfonts}
\usepackage{amscd}
\usepackage{amsxtra}
\usepackage{color}
\usepackage{enumitem}
\usepackage[
linkcolor=blue, citecolor=green, urlcolor=blue]{hyperref}

\makeatletter
\renewcommand{\@biblabel}[1]{[#1]}  
\makeatother

\newtheorem{theorem}{Theorem}[section]
\newtheorem{lemma}[theorem]{Lemma}
\newtheorem{corollary}[theorem]{Corollary}
\newtheorem{proposition}[theorem]{Proposition}

\theoremstyle{remark}
\newtheorem{remark}[theorem]{Remark}

\theoremstyle{plain}
\newtheorem{conj}{Conjecture}[section]
\theoremstyle{remark}
\newtheorem{prob}{Problem}

\theoremstyle{definition}
\newtheorem{defn}{Definition}
\newtheorem{example}{Example}

\newcommand{\C}{\mathbb C}
\newcommand{\D}{\mathbb D}

\newcommand{\Vol}{\operatorname{Vol}}
\newcommand{\id}{\operatorname{Id}}

\newcommand{\dd}{\,d}

\newenvironment{pf}[1][]{%
	\vskip 1mm
	\noindent
	\ifthenelse{\equal{#1}{}}%
	{{\slshape Proof. }}%
	{{\slshape #1.} }%
}%
{\qed\bigskip}

\newcounter{alphabet}
\newcounter{tmp}


\providecommand{\id}{{\operatorname{id}}}

\def\be{\begin{equation}}
	\def\ee{\end{equation}}

\newcommand{\bee}{\begin{enumerate}}
	\newcommand{\eee}{\end{enumerate}}

\newcommand{\blem}{\begin{lem}}
	\newcommand{\elem}{\end{lem}}
\newcommand{\bthm}{\begin{thm}}
	\newcommand{\ethm}{\end{thm}}
\newcommand{\bcor}{\begin{cor}}
	\newcommand{\ecor}{\end{cor}}
\newcommand{\beg}{\begin{example}}
	\newcommand{\eeg}{\end{example}}
\newcommand{\begs}{\begin{examples}}
	\newcommand{\eegs}{\end{examples}}
\newcommand{\bdefe}{\begin{defn}}
	\newcommand{\edefe}{\end{defn}}
\newcommand{\bprob}{\begin{prob}}
	\newcommand{\eprob}{\end{prob}}
\newcommand{\bques}{\begin{ques}}
	\newcommand{\eques}{\end{ques}}
\newcommand{\bei}{\begin{itemize}}
	\newcommand{\eei}{\end{itemize}}
\newcommand{\bcon}{\begin{conj}}
	\newcommand{\econ}{\end{conj}}
\newcommand{\bcons}{\begin{conjs}}
	\newcommand{\econs}{\end{conjs}}
\newcommand{\bprop}{\begin{propo}}
	\newcommand{\eprop}{\end{propo}}
\newcommand{\br}{\begin{rem}}
	\newcommand{\er}{\end{rem}}
\newcommand{\brs}{\begin{rems}}
	\newcommand{\ers}{\end{rems}}
\newcommand{\bo}{\begin{obser}}
	\newcommand{\eo}{\end{obser}}
\newcommand{\bos}{\begin{obsers}}
	\newcommand{\eos}{\end{obsers}}
\newcommand{\bpf}{\begin{pf}}
	\newcommand{\epf}{\end{pf}}
\newcommand{\ba}{\begin{array}}
	\newcommand{\ea}{\end{array}}
\newcommand{\beq}{\begin{eqnarray}}
	\newcommand{\beqq}{\begin{eqnarray*}}
		\newcommand{\eeq}{\end{eqnarray}}
	\newcommand{\eeqq}{\end{eqnarray*}}

\numberwithin{equation}{section}
\newcounter{minutes}
\divide\time by 60
\newcounter{hours}
\multiply\time by 60 \addtocounter{minutes}{-\time}
\begin{document}
	\bibliographystyle{plain} 
	
	
	%
	
	\title[Equality in the multidimensional Suita inequality on convex domains]{Equality in the multidimensional Suita inequality on convex domains}

	\date{\today}
	
	\author[Zhi-Gang Wang]{Zhi-Gang Wang}
	\vskip.20in
\address{\noindent Zhi-Gang Wang  \vskip.05in
School of Mathematics and Statistics, and Hunan Provincial University Key Laboratory for Big Data Analysis and Application, Hunan First Normal University,
Changsha 410205, Hunan, P. R. China.}
\email{\textcolor[rgb]{0.00,0.00,0.84}{sjyzhigangwang$@$hnfnu.edu.cn}}

	\subjclass[2020]{Primary 32A25, 32U35; Secondary 32F45, 32W05.
	}
	\keywords{Suita inequality, Bergman kernel, Azukawa
indicatrix, pluricomplex Green function, convex domain.
	}
	\begin{abstract}
We study equality conditions for the multidimensional Suita inequality on convex domains in $\mathbb C^n$. Let $D\subset\mathbb C^n$ be convex and suppose its Bergman kernel is positive at $p\in D$. Equality at $p$ holds if and only if $D$ is biholomorphic to its Azukawa indicatrix at $p$ by a map sending $p$ to the origin with identity differential at $p$. This settles the convex-domain conjecture of B\l ocki and Zwonek [New York J. Math. 21, 151--161 (2015)], covering also unbounded convex domains. Our argument extracts a holomorphic radial vector field near the Green pole via weak limits from weighted Bergman estimates, and the complex extreme points of the Azukawa indicatrix give the global biholomorphism. For smoothly bounded strongly convex domains, vanishing of the scalar Bergman curvature at one Green level set implies equality. In addition, we represent the Suita defect as an integral of this curvature and classify the admissible equality sets in the smooth setting.
\end{abstract}
	\maketitle
	\pagestyle{myheadings}
	\markboth{ Z.-G. Wang}{Equality in the multidimensional Suita inequality on convex domains}
	\section{Introduction and statements of the main results
	}	

\subsection{Background and motivation}
For a Greenian Riemann surface, Suita \cite{S1972} conjectured a fundamental relation between the logarithmic capacity $c_\beta$ and the diagonal Bergman kernel:
\[
c_\beta^2(p)\leq \pi K(p).
\]
In a local coordinate $z$ centered at $p$, this logarithmic capacity is given by
\[
\exp\bigl(\lim_{q\to p}\bigl(G(q,p)-\log|z(q)|\bigr)\bigr),
\]
where $G$ denotes the negative Green function. Here the Bergman kernel corresponds to the space of holomorphic one-forms with respect to that same local coordinate.

On Riemann surfaces, Suita's conjecture connects potential-theoretic capacity to $L^2$ holomorphic function theory.
Ohsawa \cite{O1995,OT1987} first linked the Suita inequality to the $L^2$-extension phenomenon and obtained a lower bound with a nonoptimal constant.
B\l ocki \cite{Blo2013} proved the sharp Suita inequality for planar domains, while Guan and Zhou \cite{GZ2012} settled it for arbitrary open Riemann surfaces using optimal $L^2$ extension techniques.
Subsequently, Guan and Zhou \cite{GZ2015} solved the equality problem in dimension one: equality at a single point forces the surface to be biholomorphic to the unit disk minus a relatively closed polar set.
Dong \cite{D2022} later provided an alternative proof of this rigidity result via the Maitani--Yamaguchi variation formula.
Dong, Treuer, and Zhang \cite[Theorems~1.5 and~1.6]{DTZ2023} further derived rigidity statements for planar volume inequalities on Green sublevel sets, as well as comparison relations among higher-dimensional Bergman kernels, Azukawa volumes, and Euclidean boundary distances.

While the one-dimensional equality theory rests on potential-theoretic variational tools, moving to several complex variables introduces substantial new difficulties.
There is no direct analogue of the classical Riemann-surface potential-theory framework, and compactness arguments frequently break down for unbounded domains, making pointwise-to-global rigidity far less accessible.

In several variables, let $G_D(\cdot,p)$ denote the pluricomplex Green function and $I_D^A(p)$ the Azukawa indicatrix, whose precise definition is recalled in Section~\ref{sec2}.
The pluricomplex Green function was introduced by Klimek \cite{K1985}, and Azukawa indicatrix by Azukawa \cite{A1987}.
B\l ocki \cite{Blo2014} established a sharp lower bound for $K_D(p)$ in terms of volumes of Green sublevel sets.
Using rescaled Green sublevel sets, B\l ocki and Zwonek \cite{BZ2015} obtained the multidimensional Suita inequality
\[
K_D(p)\geq \frac{1}{\operatorname{Vol}(I_D^A(p))}
\]
on pseudoconvex domains. They tested this inequality on convex ellipsoids and extended the framework to higher-order Bergman kernels; see, e.g.,\ \cite{BZ2016,BZ2020}.
Balakumar, Borah, Mahajan, and Verma \cite{BBMV2019,BBMV2020} investigated the boundary behaviour of the corresponding volume product for the Kobayashi and other indicatrices over strongly pseudoconvex, finite-type, and polyhedral model domains.
Nikolov and Thomas \cite{NT2019} related these volume comparisons to the Carath\'eodory--Eisenman volume invariants.

In 2015, B\l ocki and Zwonek \cite[Conjecture~3, p.~155]{BZ2015} formulated a natural rigidity conjecture for convex domains.

\begin{conj}
Given a convex domain $D\subset\mathbb C^n$ and a point $p\in D$ satisfying $K_D(p)>0$, the equality $K_D(p)\operatorname{Vol}(I_D^A(p))=1$ should hold if and only if there exists a biholomorphism sending $D$ onto a balanced domain while mapping $p$ to the origin.
\end{conj}
The conjecture allows the target balanced domain to be non-convex and also includes the challenging case of unbounded convex $D$, where compactness arguments are generally unavailable.

Variational properties of the Bergman kernel offer a promising route toward such equality characterizations.
Berndtsson \cite{B2006,B2009} proved plurisubharmonicity of the logarithm of the diagonal weighted Bergman kernel for pseudoconvex families and built curvature theory for the associated Hilbert vector bundles.
Berndtsson and Lempert \cite{BL2016} exploited this bundle positivity to prove sharp $L^2$-extension estimates.
Wang's curvature formula \cite[Theorem~2.17]{W2017} applies to smooth families of smoothly bounded strongly pseudoconvex domains and, under additional hypotheses, relates flatness of the full Bergman bundle to holomorphic triviality.
Guan and Mi \cite{GM2023}, Guan and Yuan \cite{GuY2023}, and Guan, Mi, and Yuan \cite{GMY2024} established concavity theorems for minimal $L^2$ integrals and studied degeneration toward linear behaviour.
Xu and Zhou \cite[Proposition~6.1 and Theorems~7.1-7.2]{XZ2026} showed that equality in the multidimensional Suita inequality enforces the minimising holomorphic function to achieve its extremal value on every Green sublevel of a bounded pseudoconvex domain; they also constructed explicit Hartogs-domain examples where equality occurs.
Very recently, Deng, Huang, and Liu \cite{DHL2026} derived period bounds for curvature of Bergman evaluation matrices on smooth planar domains. For scalar evaluation at a Green pole, equality in their bound at a regular Green level is equivalent to simple connectivity of the corresponding sublevel set.

For convex domains, the geometry of Green sublevel sets is governed by Lempert's theory of extremal holomorphic disks; see e.g.,\ \cite{L1981}.
Bracci and Saracco \cite{BS2009} characterised Kobayashi hyperbolicity for unbounded convex domains and its interplay with the Bergman metric.
Kaup and Upmeier \cite{KU1976} described the origin orbit inside the unit ball of a complex norm, while Braun, Kaup, and Upmeier \cite{BKU1978} developed the corresponding structure theory for bounded circular domains.
For smoothly bounded domains whose boundary carries a positive-definite real second fundamental form, Patrizio \cite[Theorems~5-7]{P1983} used Lempert's Monge--Amp\`ere exhaustion to give criteria for circular realisations and supplied a characterisation of the unit ball via rotations about two distinct interior points.
Further related rigidity statements rely on real and complex Monge--Amp\`ere equations and parabolic exhaustions; see e.g.,\ \cite{P1984,P1985}.

	\subsection{Main results}

Let $K_D(z,w)$ be the Bergman kernel with respect to Euclidean
volume $\dd V$, and let $G_D(\cdot,p)$ be the pluricomplex Green
function with logarithmic pole at $p$. Let $I_D^A(p)$ denote the
Azukawa indicatrix and $\mu=\mu_I$ its Minkowski functional. Set
\[
 G=G_D(\cdot,p),\quad I=I_D^A(p),\quad D_s=\{z\in D:G(z)<s\}
\]
and
\[
 F_p(s)=\log K_{D_s}(p,p)+2ns\quad(s<0).
\]
Recall that a domain $E$ is said to be \emph{balanced} if
$\lambda E\subset E$ for every $|\lambda|\le1$; see e.g., \cite{JP2013}.

We now state our main results.

\begin{theorem}\label{thm1}
Let $D\subset\mathbb C^n$ be convex, and let $p\in D$ satisfy
$K_D(p,p)>0$. Then the following statements are equivalent.
\begin{enumerate}[label=\textup{(\roman*)}]
\item $K_D(p,p)\Vol(I)=1$;
\item There exist $s_1<s_2<0$ with $F_p(s_1)=F_p(s_2)$;
\item There is a biholomorphism $\Psi:I\to D$ such that
      $\Psi(0)=p$, $d\Psi_0=\id$, and
      $G(\Psi(w))=\log\mu_I(w)$;
\item There is a biholomorphism from $D$ onto a balanced domain
      taking $p$ to the origin.
\end{enumerate}
If these conditions fail, $F_p$ is strictly increasing on $(-\infty,0)$.
\end{theorem}

The equivalence of (i) and (iv) answers the conjecture posed by B\l ocki and Zwonek \cite[Conjecture~3, p.~155]{BZ2015}. The fixed domains in the
weighted estimates are relatively compact Green sublevels. This
allows the proof to include unbounded $D$.

We call a bounded convex domain with $C^\infty$ boundary
\emph{strongly convex} if the real Hessian of a smooth defining
function is positive definite on each real tangent space of
the boundary. In this case, the Lempert circular representation
\[\Phi_p:I\to D\] is defined by
\[\Phi_p(\zeta v)=\phi_v(\zeta)\] for $v\in\partial I$, where
$\phi_v$ is the unique complex geodesic with
$\phi_v(0)=p$ and $\phi_v'(0)=v$. This representation is a
homeomorphism and is smooth off the origin; see e.g., \cite{L1981,P1983}.

\begin{theorem}\label{thm2}
Let $D\subset\mathbb C^n$ be bounded, $C^\infty$-smooth and strongly
convex, and let $p\in D$. Then the following statements are equivalent.
\begin{enumerate}[label=\textup{(\roman*)}]
\item $K_D(p,p)\Vol(I)=1$;
\item $F_p''(s_0)=0$ for some $s_0<0$;
\item $\Phi_p:I\to D$ is biholomorphic;
\end{enumerate}
If these conditions do not hold, $F_p''(s)>0$ and $F_p'(s)>0$ for every $s<0$.
Moreover,
\begin{equation}\label{eq:defect}
 \log\!\bigl(K_D(p,p)\Vol(I)\bigr)
   =\int_{-\infty}^{0}(-s)F_p''(s)\dd s .
\end{equation}
\end{theorem}

The proof of Theorem~\ref{thm1} starts with the common normalized
Bergman extremal on the Green sublevels. A quadratic cutoff
in the Green variable gives a weighted kernel on a fixed domain
whose limit is constant in the level parameter. The second variation
controls the antiholomorphic derivatives of an adjoint solution
and its error in the radial identity. Theorem~\ref{thm2} uses the
moving-domain formula of Wang \cite{W2017} on a single fiber.

\section{Green functions and the common extremal}\label{sec2}

For a domain $D\subset\C^n$, write $A^2(D)$ for the Hilbert space of
holomorphic functions square integrable with respect to $\dd V$.
The diagonal Bergman kernel is characterized by
\[
 K_D(p,p)=\sup\{|f(p)|^2:f\in A^2(D),\ \|f\|_2\le1\}.
\]
When it is positive, its reciprocal is the minimum of $\|f\|_2^2$
under $f(p)=1$. The unique minimizer is
$K_D(\cdot,p)/K_D(p,p)$. For a positive weight, we replace $\dd V$
by the weighted volume; see e.g., \cite{B2006,BZ2015}.

The pluricomplex Green function $G_D(z,p)$ is the upper envelope of
negative plurisubharmonic functions $u$ on $D$ satisfying
\[u(z)\le\log|z-p|+O(1)\] near $p$; see e.g., \cite{K1991,BZ2015}. The Azukawa
indicatrix is
\[
 I_D^A(p)=\left\{v\in\C^n:
  \limsup_{\lambda\to0}\bigl(G_D(p+\lambda v,p)
                            -\log|\lambda|\bigr)<0\right\}.
\]
The Azukawa pseudometric is obtained by taking the exponential of this limsup-quantity; see e.g., 
\cite{A1987,BZ2015}. For an open balanced set $I$ containing $0$,
its Minkowski functional is \[\mu_I(v)=\inf\{r>0:v\in rI\}.\] If $I$
is bounded and convex, $\mu_I$ is a norm over $\C$.

The Carath\'eodory distance $c_D$ is the supremum of Poincar\'e distances on the unit disk $\mathbb{D}$ pulled back by holomorphic maps $f\colon D\to\mathbb{D}$. We adopt the standard normalization
\[
c_{\mathbb{D}}(0,r) = \operatorname{arctanh} r
\]
of the Poincar\'e distance. Related rigidity results for balanced domains in this setting may be found in~\cite{Wang26}.
The Lempert
function is the infimum of Poincar\'e distances of pairs on a single
holomorphic disk in $D$. The Kobayashi distance $k_D$ is the largest
pseudodistance bounded above by the Lempert function. For $v\ne0$,
the infinitesimal quantities $\gamma_D(p;v)$ and $\kappa_D(p;v)$ are
the supremum of $|dh_p(v)|$ over $h:D\to\D$ with $h(p)=0$, and the
infimum of $1/r$ over disks $\phi$ with $\phi(0)=p$ and
$\phi'(0)=rv$ ($r>0$), respectively. A \emph{complex geodesic} here
is a holomorphic disk $\phi$ with a holomorphic left inverse
$\ell:D\to\D$; see e.g., \cite{JP2013,L1981}.

\begin{lemma}\label{lem1}
Under the hypotheses of Theorem~\ref{thm1}, $I$ is a bounded convex
balanced domain, and
\begin{equation}\label{eq:green-metric}
 e^{G(z)}=\tanh c_D(p,z)
 =\sup\{|h(z)|:h\in\mathcal O(D,\D),\ h(p)=0\}.
\end{equation}
In a neighborhood of $p$ one has
\begin{equation}\label{eq:green-tangent}
 e^{G(p+w)}=\mu_I(w)+O(|w|^2)
\end{equation}
uniformly in $w$. The function $e^G$ is locally Lipschitz, with
value zero at $p$. Each negative Green sublevel $D_s$ is convex
and relatively compact, and
\begin{equation}\label{eq:sublevel}
 G_{D_s}(z,p)=G(z)-s,\quad I_{D_s}^A(p)=e^s I.
\end{equation}
For every $v\in\partial I$, there is a complex geodesic $\phi_v$
with $\phi_v(0)=p$ and $\phi_v'(0)=v$.
\end{lemma}

\begin{proof}
By the decomposition of Bracci and Saracco \cite[Proposition~1.2]{BS2009},
a convex domain containing a complex affine line is linearly
equivalent to a product with a positive-dimensional complex
Euclidean factor. Fubini's theorem and the mean-value inequality
give $A^2(D)=\{0\}$ in this case. Positivity of the kernel therefore
excludes complex lines. Bracci and Saracco's theorem gives tautness,
compact closed Kobayashi balls, and convex open Kobayashi balls
\cite[Propositions~3.2 and~3.4]{BS2009}. Lempert's theorem identifies the Lempert function, $k_D$ and $c_D$; see e.g., \cite{BS2009,L1981}. It also gives $\gamma_D=\kappa_D$. In the
unbounded case, the infinitesimal equality follows by exhausting
with $D\cap B(p,R)$. Normalized disk-valued extremals and Montel's
theorem give convergence on the Carath\'eodory side. On the
Kobayashi side one tests disks restricted to $r\D$ and lets
$r\uparrow1$.

The functions $\log|h|$ are competitors for $G$. For the opposite
bound, pull a competitor back along a disk through $p$ and apply
the subharmonic Schwarz inequality with a logarithmic pole.
Lempert's distance equality gives \eqref{eq:green-metric}. Uniform
Cauchy estimates give local Lipschitz regularity of the supremum
there and a first-order expansion at $p$ with leading term
$\gamma_D(p;w)$. In particular, we have \[\gamma_D(p;w)=\mu_I(w).\]
This is a norm equivalent to the Euclidean norm, since a convex
domain without complex lines admits $n$ independent complex
linear support functionals into half-planes \cite{BS2009}.
Tautness gives extremal disks for the point and tangent problems by
normal limits. The normalized Carath\'eodory functions attain their
extremal values by Montel's theorem. Schwarz equality makes
them left inverses of the corresponding disks.

For $z\in D_s\setminus\{p\}$, choose a two-point extremal disk
$\phi$ through $p,z$ with left inverse $\ell$, normalized at $p$.
By \eqref{eq:green-metric}, we have \[G(\phi(\zeta))=\log|\zeta|.\]
The function $e^{-s}\ell$ is disk-valued on $D_s$. The Green lower
bound from this function and the upper bound from the disk
$\phi(e^s\,\cdot)$ give the first identity in \eqref{eq:sublevel}.
The second follows from the definition of the indicatrix. The
negative Green sublevels are Kobayashi balls by
\eqref{eq:green-metric}.
\end{proof}

\begin{proposition}\label{prop1}
Put $V=\Vol(I)$. Under the hypotheses of Theorem~\ref{thm1}, $F_p$
is convex and nondecreasing, with
\begin{equation}\label{eq:limits}
 F_p(-\infty)=-\log V,\quad F_p(0^-)=\log K_D(p,p).
\end{equation}
If either (i) or (ii) of Theorem~\ref{thm1} holds, there exist
$a<b<0$ and $f\in A^2(D_b)$, $f(p)=1$, whose restriction is
the normalized Bergman extremal on each $D_u$, $u\le b$, and
\begin{equation}\label{eq:common}
 \int_{D_u}|f|^2\dd V=V e^{2nu}\quad(u\le b).
\end{equation}
\end{proposition}

\begin{proof}
The set \[\{(t,z)\in\C\times D:{\rm Re} t<0,\ G(z)<{\rm Re} t\}\] is
pseudoconvex. Indeed, it admits an exhaustion obtained by adding
a nonnegative plurisubharmonic exhaustion of $D$ to
\[|t|^2-\log(-{\rm Re} t)-\log(1-e^{G-{\rm Re} t}),\] using a constant lower
bound for the first two terms in $t$. The last term is a convex
increasing function of $G-{\rm Re} t$. Berndtsson's theorem
\cite[Theorem~1.1]{B2006} gives convexity of
\[s\mapsto\log K_{D_s}(p,p)\] by taking the logarithmic diagonal
kernel at the constant section $p$ and using independence of ${\rm Im} t$.

By \eqref{eq:green-tangent}, $e^{-s}(D_s-p)$ lies between
$(1-Ce^s)I$ and $(1+Ce^s)I$ for sufficiently negative $s$.
Indeed, the nested compact sublevels shrink to $p$ and the
Euclidean and indicatrix norms are equivalent. Inclusion
monotonicity and the dilation law give the first limit in
\eqref{eq:limits}, since $K_I(0,0)=1/V$ by rotation averaging.
The other limit follows from increasing exhaustion, using
restriction of extremals for one inequality and a normal-family
and Fatou argument for the other. Convexity and the finite limit
at $-\infty$ give monotonicity. If the endpoint values agree,
$F_p$ is constant. If $F_p$ has equal values at two distinct finite
levels, convexity makes it constant up to the larger level.

Take $a<b<0$ inside the resulting constant interval. By
\eqref{eq:sublevel}, $D_b$ attains the Suita equality at $p$.
Apply the degree-zero case of Xu and Zhou's equality result
\cite[Proposition~6.1(2)]{XZ2026} to the bounded pseudoconvex
domain $D_b$. It supplies a single normalized extremal on all
its Green sublevels, which by \eqref{eq:sublevel} are the
$D_u$ for $u<b$. The constant value $-\log V$ gives
\eqref{eq:common}.
\end{proof}

\begin{remark}\label{rem1}
The defect also has a measure formula on arbitrary convex domains
under the hypotheses of Theorem~\ref{thm1}. Let $\nu_p$ be the positive
distributional second derivative of $F_p$ on $(-\infty,0)$. The
finite left limit of $F_p$ gives \[F'_{p,+}(-\infty)=0.\] Here
$F'_{p,+}$ is the right derivative of the convex function.
Consequently, we have \[F'_{p,+}(s)=\nu_p((-\infty,s]).\] Local absolute
continuity of $F_p$ and Tonelli's theorem give
\[
 \log\!\bigl(K_D(p,p)\Vol(I)\bigr)
   =\int_{(-\infty,0)}(-s)\,\dd\nu_p(s).
\]
In particular, the two-level condition in Theorem~\ref{thm1} detects
vanishing of this measure on a left interval. The single-level
criterion in Theorem~\ref{thm2} uses the smooth boundary geometry
of the moving fibers.
\end{remark}

\section{Weighted kernel variation}\label{sec3}

We prove that the equality assumption yields a holomorphic vector field near
$p$ satisfying the radial equation almost everywhere. The argument keeps the
underlying domain fixed and varies its weight.

\begin{lemma}\label{lem2}
Let $\Omega_0\Subset\mathbb C^n$ be a smoothly bounded pseudoconvex domain
containing $p$. Let $g$ be smooth and plurisubharmonic on a neighborhood
of $\overline{\Omega_0}$, and let $\chi$ be smooth, convex, and nondecreasing.
For real $s$ consider the weight
\[
 \psi_s(z)=\chi(g(z)-s)+2ns.
\]
Write $P_s$ for the weighted Bergman projection, $u_s=K_s(\cdot,p)$ for the
weighted kernel section, $k_s=K_s(p,p)$, and $\|\cdot\|_s$ for the weighted
$L^2$ norm with measure $e^{-\psi_s}\dd V$. There is a $(0,1)$-form
\[\alpha_s=\sum_a\alpha_{s,a}\,\dd\bar z_a\] such that
\begin{equation}\label{eq:uniform-adjoint}
 \|\alpha_s\|_s\le C\|\dot\psi_s u_s\|_s,
 \end{equation} and
 \begin{equation}\label{eq:radial-curvature-bound}
 k_s(\log k_s)''\ge
 \int_{\Omega_0}\left(
 \sum_{a,b}\left|\frac{\partial\alpha_{s,a}}{\partial\bar z_b}\right|^2
 +\chi''(g-s)\left|u_s+\sum_a g_{z_a}\alpha_{s,a}\right|^2
 \right)e^{-\psi_s}\dd V.
 \end{equation} 
The constant $C$ can be chosen to depend only on a fixed ball containing
$\Omega_0$.
\end{lemma}

\begin{proof}
We take the inner products to be linear in the first variable. On compact
parameter intervals, the weights are smooth, bounded, and uniformly bounded
away from zero on the fixed domain. Thus the weighted inner products on the
fixed unweighted Bergman space are represented by a smoothly varying family
of bounded invertible operators. In particular, the reproducing identity
can be differentiated with respect to $s$. Differentiating it once gives
\[\dot u_s=P_s(\dot\psi_su_s),\] whence
\[
 \dot k_s=\int_{\Omega_0}\dot\psi_s|u_s|^2 e^{-\psi_s}\dd V.
\]
Set \[w_s=(1-P_s)(\dot\psi_su_s).\] Differentiating the last integral gives
\[
 \ddot k_s=\int_{\Omega_0}\ddot\psi_s|u_s|^2e^{-\psi_s}\dd V
       +\|P_s(\dot\psi_su_s)\|_s^2-\|w_s\|_s^2.
\]
The Cauchy--Schwarz inequality bounds $\dot k_s^2/k_s$ by
$\|P_s(\dot\psi_su_s)\|_s^2$. Hence
\begin{equation}\label{eq:weighted-scalar-curvature}
 k_s(\log k_s)''\ge
 \int_{\Omega_0}\ddot\psi_s|u_s|^2e^{-\psi_s}\dd V-\|w_s\|_s^2.
\end{equation}
This argument uses the weighted Bergman variation framework of \cite{B2006}.

For each $\bar\partial$-closed weighted square-integrable $(0,1)$-form
$\eta$, the H\"ormander estimate with weight $\psi_s+|z|^2$ gives a solution
of $\bar\partial h=\eta$ satisfying \[\|h\|_s\le C\|\eta\|_s.\] Indeed,
the perturbed weight has complex Hessian bounded below by the Euclidean
metric, and its norm is equivalent to $\|\cdot\|_s$ with constants depending
only on the containing ball. Since $w_s$ is orthogonal to holomorphic
functions, the functional $\eta\mapsto(h,w_s)_s$ is well-defined and bounded
on the closed space $\ker\bar\partial$ of $(0,1)$-forms. Its Riesz
representative $\alpha_s$ belongs to $\ker\bar\partial$ and satisfies
\[\bar\partial^*_{\psi_s}\alpha_s=w_s\] and \eqref{eq:uniform-adjoint}.
This construction uses the weighted $L^2$ estimate in \cite{H1965}.
For a fixed smooth family one can choose $\alpha_s$ measurably in $s$.
For example, fix an unweighted bounded right inverse on $\ker\bar\partial$
and use it for $h$ in the defining functional. The weighted inner products
on this fixed closed subspace are smoothly varying positive invertible
operators, so their Riesz representatives depend continuously on $s$.

The Morrey--Kohn--H\"ormander estimate for forms in the domains of
$\bar\partial$ and $\bar\partial^*_{\psi_s}$ gives
\begin{equation}\label{eq:weighted-morrey}
 \|w_s\|_s^2\ge
 \int_{\Omega_0}\left(
 \sum_{a,b}|\partial_{\bar z_b}\alpha_{s,a}|^2+
 \sum_{a,b}(\psi_s)_{z_a\bar z_b}\alpha_{s,a}\overline{\alpha_{s,b}}
 \right)e^{-\psi_s}\dd V.
\end{equation}
The boundary term is nonnegative by pseudoconvexity. In particular, the
first integrand is the full array of coefficient derivatives, even though
the exterior derivative $\bar\partial\alpha_s$ vanishes. The inequality
on the indicated operator domains follows from the smooth identity via
graph-norm approximation; see \cite{ST2010}.

By the adjoint equation
\[\|w_s\|_s^2=(\bar\partial w_s,\alpha_s)_s\]
and \[\bar\partial w_s=u_s\bar\partial\dot\psi_s,\] we can
write \[-\|w_s\|_s^2=\|w_s\|_s^2-2{\rm Re}(u_s\bar\partial\dot\psi_s,\alpha_s)_s\]
in \eqref{eq:weighted-scalar-curvature}. Apply \eqref{eq:weighted-morrey}
to the positive term on the right. The Hessian of the weight in $z$ is
\[\chi''(g-s)\partial g\otimes\bar\partial g+\chi'(g-s)\partial\bar\partial g,\]
whereas \[\ddot\psi_s=\chi''(g-s)\] and
\[\bar\partial\dot\psi_s=-\chi''(g-s)\bar\partial g.\] The $\chi''$ terms
combine into the square in \eqref{eq:radial-curvature-bound}. The remaining
Hessian term with coefficient $\chi'$ is nonnegative and can be discarded.
\end{proof}

\begin{lemma}\label{lem3}
Suppose that the Green function $G$ is locally Lipschitz away from $p$,
its negative sublevels are relatively compact and pseudoconvex, and
there is an $f$ on some $D_b$ which has the common-extremal property
in Proposition~\ref{prop1} and satisfies \eqref{eq:common}. Then there is a
holomorphic vector field $X$ on a neighborhood of $p$ such that
\begin{equation}\label{eq:weak-radial-field}
 \partial G(X)=\frac12
\end{equation} almost everywhere away from $p$.
\end{lemma}

\begin{proof}
Take $s_*<b$ from the range in the common-extremal assumption
and set $\Omega=D_{s_*}$. Fix a relatively compact open real
parameter interval $J\Subset(-\infty,s_*)$. Let $j\ge2$. Define a convex nondecreasing function
$\widetilde\chi_j$ by setting it equal to $-2\log(-u)$ on
$[-j,-1/j]$, extending as a constant on $u\le-j$ and by the tangent line
at $-1/j$ on $u\ge-1/j$. Convolve with a nonnegative smooth mollifier of
radius $j^{-3}$. Denote the resulting smooth function by $\chi_j$, and set
$q_j=e^{-\chi_j}$. These functions satisfy
\begin{equation}\label{eq:weight-cutoff-bounds}
 q_j(u)\longrightarrow (-u)_+^2, \quad
 0\le q_j(u)\le C(1+(-u)_+^2),\quad
 |\chi_j'(u)|^2q_j(u)\le C\end{equation}
 and
 \begin{equation}\label{11}
 \chi_j''(u)q_j(u)\longrightarrow2
 \ \ \text{locally uniformly on }(-\infty,0).
\end{equation}
Indeed, the derivative-squared estimate is immediate on the logarithmic
part and on each extension. In the mollification neighborhood of either
join, the adjacent values of $-u$ are comparable to $j$ or $1/j$,
respectively, and the same estimate follows from the adjacent derivative
bounds. The locally uniform limit follows by differentiating the smooth
logarithmic profile on compact subintervals. The tangent extension also
shows that $q_j(u)\le Cj^{-2}$ for $u\ge0$.

Consider first the unsmoothed weight
\[\psi_{j,s}=\chi_j(G-s)+2ns\] on $\Omega$. By \eqref{eq:common}, the
function $f$ with $f(p)=1$, minimizes the squared norm on every Green
sublevel and on $\Omega$ itself. Since $q_j$ is nonincreasing, layer-cake
integration expresses any norm for the weight $\psi_{j,s}$ as a
nonnegative combination of these unweighted squared norms. Explicitly,
\[
 q_j(G-s)=q_j(s_*-s)+\int_{-\infty}^{s_*-s}
                    [-q_j'(a)]\mathbf1_{\{G<s+a\}}\,\dd a
\] almost everywhere on $\Omega$.
Thus $f$ is the weighted extremal for every $(j,s)$. Write $k_{j,s}$ for
the weighted diagonal kernel on $\Omega$. Its kernel section is $k_{j,s}f$.
The measure $|f|^2\dd V$ on $\Omega$, pushed forward by $G$, has density
$2nV e^{2n\sigma}\dd\sigma$ on $(-\infty,s_*)$ by \eqref{eq:common}.
Consequently, we have
\begin{equation}\label{eq:weighted-common-kernel}
 k_{j,s}^{-1}=\|f\|_{\psi_{j,s}}^2
 =2nV\int_{-\infty}^{s_*-s}q_j(a)e^{2na}\,\dd a
 \longrightarrow \frac{V}{2n^2},
\end{equation}
uniformly on compact subintervals of $J$. To see the uniformity, all upper
limits on such an interval are in a compact subset of $(0,\infty)$, and
the positive tail of the integral is uniformly $O(j^{-2})$. The part up
to zero converges by dominated convergence. In addition,
\begin{equation}\label{eq:weighted-common-derivative}
 \sup_{j\ge2,\ s\in J}\int_{\Omega}|2n-\chi_j'(G-s)|^2
 |k_{j,s}f|^2 e^{-\psi_{j,s}}\dd V<\infty.
\end{equation}
This estimate is a consequence of \eqref{eq:weight-cutoff-bounds}, \eqref{11}, \eqref{eq:weighted-common-kernel},
together with the assumption $f\in L^2(\Omega)$.

Let $\Omega_\ell\Subset\Omega$ be a smooth bounded pseudoconvex exhaustion
and let $G_\ell$ be obtained by convolution of $G$ on a neighborhood of
$\overline\Omega$ with a nonnegative mollifier of radius tending to zero.
Then $G_\ell$ is smooth and plurisubharmonic, it tends to $G$ almost
everywhere and locally uniformly off $p$, and its gradient tends to that
of $G$ in $L^2$ on compact subsets off $p$. For fixed $j$, apply
Lemma~\ref{lem2} on $\Omega_\ell$ with $g=G_\ell$ and weight
\[\psi_{j,\ell,s}=\chi_j(G_\ell-s)+2ns.\] Denote the corresponding kernel
section, diagonal, and adjoint form by $u_{j,\ell,s}$, $k_{j,\ell,s}$,
and $\alpha_{j,\ell,s}$.

For clarity, we justify the approximation of the weighted kernels. For
fixed $j$, all these weights are uniformly bounded above and below on
$\Omega_\ell$ and on compact parameter intervals, because $G_\ell$ is
uniformly bounded above and $\chi_j$ has a constant left tail. For fixed
$s$, restriction of the exact normalized weighted extremal bounds the
limsup of the minimum squared norm on $\Omega_\ell$ by the minimum on
$\Omega$. Conversely, the uniform local bounds on the normalized
approximating extremals give a normal-family limit on $\Omega$, and
Fatou's lemma bounds its squared norm by the liminf of the minima. It
follows that $k_{j,\ell,s}\to k_{j,s}$ and $u_{j,\ell,s}\to k_{j,s}f$
locally uniformly. The square-root weighted sections, extended by zero
from $\Omega_\ell$, converge strongly in $L^2(\Omega)$, since they
converge almost everywhere, are bounded in $L^2$, and their squared norms
are their diagonals, which converge. As $\dot\psi_{j,\ell,s}$ is uniformly
bounded for fixed $j$, the boundedness and pointwise convergence of this multiplier give
convergence of $\|\dot\psi_{j,\ell,s}u_{j,\ell,s}\|^2_{\psi_{j,\ell,s}}$
to its unsmoothed counterpart. For fixed $j$ these conclusions can be
integrated in $s$ by dominated convergence. A common interior ball
about $p$ and the uniform upper and lower weight bounds provide upper
and lower bounds for the diagonals on compact parameter intervals.

Choose an increasing sequence $\ell(j)$ such that the preceding integrated
convergences have errors tending to zero on compact parameter intervals,
and such that convergence of the kernel sections in local joint $L^2$
and convergence of the gradients of $G_{\ell(j)}$ on compact annuli
hold along the same sequence. For the latter convergence one may also
simply require $\ell(j)\to\infty$. Suppress the $\ell(j)$ index and
write $u_{j,s}, k^*_{j,s},\alpha_{j,s},\psi^*_{j,s},g_j$ for the smooth
objects. By \eqref{eq:weighted-common-kernel}, we have
\[\log k^*_{j,s}\to\log(2n^2/V)\] locally in $L^1(J)$.
The functions $\log k^*_{j,s}$ are convex by
\eqref{eq:radial-curvature-bound}, so they are locally uniformly bounded.
For a nonnegative $\eta\in C_c^\infty(J)$, this gives
\begin{equation}\label{eq:curvature-limit-zero}
 0\le\int_J\eta(s)k^*_{j,s}(\log k^*_{j,s})''\,\dd s
 \le C_\eta\int_J\eta''(s)\log k^*_{j,s}\,\dd s\longrightarrow0.
\end{equation}
In the middle inequality, use the local upper bound for $k^*_{j,s}$ and
the nonnegativity of the second derivative.

Choose $J_0\Subset J$ of positive length and a ball $B_0$ centered at $p$
with closure in $\Omega$, so small that $g_j\le\inf J_0-\delta$ on $B_0$
for all sufficiently large $j$ and some $\delta>0$. Such a choice is
possible because $G$ tends to $-\infty$ at $p$ and each $g_j$ is a local
average of $G$. The weights $e^{-\psi^*_{j,s}}$ have a positive lower bound
there, independent of $j$ and $s\in J_0$, by monotonicity of $\chi_j$
and its convergence at the fixed point $-\delta$. The bound
\eqref{eq:uniform-adjoint}, the integrated convergence of derivative
norms, and \eqref{eq:weighted-common-derivative} give a subsequence
$\alpha_{j,s}\rightharpoonup\alpha_s$ in joint $L^2(J_0\times B_0)$.
The first integrand of \eqref{eq:radial-curvature-bound} and
\eqref{eq:curvature-limit-zero} give
$\partial_{\bar z_b}\alpha_{s,a}=0$ in distributions on this product.

On each compact annulus in $B_0\setminus\{p\}$, the quantities $g_j-s$
range in a compact subset of $(-\infty,0)$, so the factor
\[\chi_j''(g_j-s)e^{-\psi^*_{j,s}}\] is uniformly bounded away from zero
by \eqref{eq:weight-cutoff-bounds}. The second integrand of
\eqref{eq:radial-curvature-bound} thus gives
\[u_{j,s}+\sum_a(g_j)_{z_a}\alpha_{j,s,a}\to0\] in $L^2$ on each such
parameter--annulus product. The kernel sections converge in \(L^2\) to
\[u = \frac{2n^2}{V} f,\]
and the product with the gradient converges in the sense of distributions. This follows from the strong \(L^2\)-convergence of the gradients together with the weak \(L^2\)-convergence of the forms. We then obtain
\[u+\sum_{a} G_{z_a}\alpha_{s,a}=0\]
almost everywhere on \(J_0\times (B_0\setminus\{p\})\).

By testing against a countable dense family and using Fubini, choose one
$s\in J_0$ for which the coefficient functions of $\alpha_s$ are weakly
holomorphic on $B_0$ and the last identity holds almost everywhere off
$p$. The weakly holomorphic coefficients have holomorphic representatives.
After shrinking $B_0$, $u$ is nonvanishing since $f(p)=1$. The holomorphic
vector field
\[
 X=-\frac1{2u}\sum_{a=1}^n\alpha_{s,a}\frac{\partial}{\partial z_a}
\]
then satisfies \eqref{eq:weak-radial-field}.
\end{proof}

\section{Local linearization and global extension}\label{sec4}

Put \(q=e^G\), with \(q(p)=0\), and write \(\mu=\mu_I\). We use the consequence of
\eqref{eq:green-metric} that
\[
 q(z)=\sup\big\{|h(z)|:h\in\mathcal O(D,\mathbb D),\ h(p)=0\big\}.
\]
In particular, Cauchy's estimate on balls compactly contained in \(D\)
makes \(q\) locally Lipschitz, including at \(p\).

\begin{lemma}\label{lem4}
Suppose that a holomorphic vector field \(X\), defined near \(p\), satisfies
\(\partial G\cdot X=1/2\) almost everywhere off \(p\). Then there is a
biholomorphic map \(H\) between neighborhoods of \(p\) and \(0\) such that
\[
 H(p)=0,\quad DH(p)=\mathrm{Id},\quad q=\mu\circ H.
\]
\end{lemma}

\begin{proof}
Let \(\phi_\tau\) be the local holomorphic flow of \(X\). The chain rule for
locally Lipschitz functions gives
\[
 \frac{d}{dt}q(\phi_t(z))=q(\phi_t(z)),\quad
 \frac{d}{dt}q(\phi_{it}(z))=0
\]
for almost every \((t,z)\) in a product neighborhood where the flows are
defined. Indeed, for each fixed \(t\) the flow maps are local
diffeomorphisms, so the exceptional set for the identity on \(\partial G\)
has null pullback. At points of differentiability off \(p\) the two
derivatives are \(2\operatorname{Re}(\partial q\cdot X)=q\) and
\(2\operatorname{Re}(i\partial q\cdot X)=0\). Fubini's theorem and absolute
continuity give the integrated identities for almost every starting point.
Continuity then gives them for every starting point. The flow property
now yields
\begin{equation}\label{eq:local-flow-q}
 q(\phi_\tau(z))=e^{\operatorname{Re}\tau}q(z)
\end{equation}
for small complex \(\tau\). Since \(q^{-1}(0)=\{p\}\), the flow fixes
\(p\), and \(X(p)=0\).

Write \(A=DX(p)\). Differentiation of the flow at its fixed point gives
\[D\phi_\tau(p)=e^{\tau A}.\] Apply \eqref{eq:local-flow-q} at \(p+rw\),
divide by \(r>0\), and use \eqref{eq:green-tangent} as \(r\downarrow0\).
This gives
\[
 \mu(e^{\tau A}w)=e^{\operatorname{Re}\tau}\mu(w).
\]
Consequently, \(e^{\tau(A-\mathrm{Id})}\) is a linear isometry of \(\mu\)
for small \(\tau\), and the exponential law extends this to all complex
\(\tau\). The norm \(\mu\) is equivalent to the Euclidean norm. Thus the
entire matrix-valued function \(e^{\tau(A-\mathrm{Id})}\) is bounded.
Liouville's theorem gives \(A=\mathrm{Id}\).

Write \[X(p+w)=w+R(w)\] with \(R(w)=O(|w|^2)\). Fix
\(0<\varepsilon<1/2\) and choose a ball about \(p\) sufficiently small that
\(|R(w)|\le\varepsilon|w|\) on its closure. Along the negative real flow, we have
\[
 \frac{d}{dt}|\phi_{-t}(z)-p|^2
 \le -2(1-\varepsilon)|\phi_{-t}(z)-p|^2.
\]
Thus the negative flow stays in the ball for all \(t\ge0\) and
\[
 |\phi_{-t}(z)-p|\le e^{-(1-\varepsilon)t}|z-p|.
\]
By continuation, \eqref{eq:local-flow-q} holds on these trajectories for
all negative real times. Furthermore,
\[
 \frac{d}{dt}\big[e^t(\phi_{-t}(z)-p)\big]
 =-e^t R(\phi_{-t}(z)-p).
\]
The right side is locally uniformly integrable in \(t\) by the decay
estimate. Hence
\[
 H(z)=\lim_{t\to\infty}e^t(\phi_{-t}(z)-p)
\]
is holomorphic, with \(H(p)=0\) and \(DH(p)=\mathrm{Id}\). It is
biholomorphic after restriction to a neighborhood of \(p\). Finally,
\eqref{eq:green-tangent} and \eqref{eq:local-flow-q} give
\[
 q(z)=e^t q(\phi_{-t}(z))
     =\mu\big(e^t(\phi_{-t}(z)-p)\big)
       +O\big(e^t|\phi_{-t}(z)-p|^2\big).
\]
The error tends to zero. This proves \(q(z)=\mu(H(z))\).
\end{proof}

\begin{proposition}\label{prop2}
Assume the hypotheses of Lemma~\ref{lem1}. If \(H\) has the
properties of Lemma~\ref{lem4}, then its local inverse extends to a
biholomorphism \(\Psi:I\to D\) satisfying
\[
 \Psi(0)=p,\quad D\Psi(0)=\mathrm{Id},\quad q(\Psi(w))=\mu(w).
\]
\end{proposition}

\begin{proof}
Let \(E\) be the set of complex extreme points of the compact convex set
\(\overline I\). In the terminology of \cite{G1975}, \(v\) is complex extreme if
\(v+\zeta a\in\overline I\) for every \(|\zeta|\le1\) implies \(a=0\).
In particular, we have \(E\subset\partial I\).

We first record two elementary properties of \(E\). Any holomorphic germ
\(\eta\) at \(0\) with values in \(\overline I\) and \(\eta(0)=v\in E\)
is constant. In fact, write
\[\eta(\zeta)=v+\sum_{m\ge1}a_m\zeta^m,\] and fix a small radius \(r>0\).
For \(|c|\le1/2\), the nonnegative density
\(1+c e^{-im\theta}+\bar c e^{im\theta}\) has normalized integral one.
Taking the corresponding average of \(\eta(re^{i\theta})\) shows that
\(v+c r^m a_m\in\overline I\). Complex extremality gives \(a_m=0\) for
each \(m\).

The second property is that for any scalar holomorphic polynomial \(P\),
\[
 \max_{\overline I}|P|=\sup_E|P|,\quad
 \max_{\overline I}\operatorname{Re}P=\sup_E\operatorname{Re}P.
\]
For either maximum choose a maximizing point of largest Euclidean norm.
If a nonconstant affine complex disk in \(\overline I\) were centered
there, the maximum principle on that disk would make the quantity being
maximized constant along it. For every sufficiently small nonzero
displacement along the disk, one of the two opposite points has strictly
larger Euclidean norm than the center. This contradicts the choice of the
maximizer, which therefore belongs to \(E\).

Let \(\Psi_0=H^{-1}\) near zero and write
\[
 \Psi_0(w)=p+\sum_{m\ge1}P_m(w),
\]
where \(P_m\) is a vector-valued homogeneous polynomial of degree \(m\).
For \(v\in E\), choose a Lempert disk \(\varphi_v:\mathbb D\to D\)
with \(\varphi_v(0)=p\) and \(\varphi_v'(0)=v\), as in
Lemma~\ref{lem1}. Its holomorphic left inverse gives
\(q(\varphi_v(\zeta))=|\zeta|\). For small \(\zeta\ne0\), the map
\[
 \zeta\longmapsto\frac{H(\varphi_v(\zeta))}{\zeta}
\]
has \(\mu\)-norm one and extends holomorphically across zero with value
\(v\). The first property of \(E\) proves that it is constant, so
\(\Psi_0(\zeta v)=\varphi_v(\zeta)\) near zero.

For \(0<r<1\), put \[C_r=\{z\in D:q(z)\le r\}.\] By
Lemma~\ref{lem1}, \(C_r\) is compact in \(D\). Applying Cauchy's coefficient
estimate to \(\varphi_v-p\) on \(|\zeta|=r\) gives
\[
 |P_m(v)|\le M_r r^{-m}\ \ (v\in E),\quad
 M_r=\max_{z\in C_r}|z-p|.
\]
The same bound holds on \(\overline I\) by the maximum property of \(E\),
applied to scalar complex linear functionals composed with \(P_m\).
For \(0<\rho<r\), homogeneity gives \[|P_m(w)|\le
M_r(\rho/r)^m\] on \(\rho\overline I\). Thus the series defines a
holomorphic extension \(\Psi\) on \(I\). The identity theorem gives
\(\Psi(\zeta v)=\varphi_v(\zeta)\) for \(v\in E\) and \(|\zeta|<1\).

For each real linear functional \(L\) on \(\mathbb C^n\), apply the
maximum property to \(L\) composed with the Taylor polynomials of
\(w\mapsto\Psi(rw)\). This is legitimate since \(L\) is the real part of
a complex linear functional. Taking the uniform limit on \(\overline I\)
gives
\[
 \sup_{w\in r\overline I} L(\Psi(w))
 =\sup_{v\in E} L(\Psi(rv))
 \le \max_{z\in C_r} L(z).
\]
Convex separation gives \[\Psi(r\overline I)\subset\operatorname{conv}(C_r).\]
The latter set is compact and is contained in \(D\), by finite-dimensional
convexity. Thus \(\Psi(I)\subset D\).

Now let \(v\) be any point of \(\partial I\). For every
\(h\in\mathcal O(D,\mathbb D)\) with \(h(p)=0\), Schwarz's lemma gives
\(|h(\Psi(\zeta v))|\le|\zeta|\). Taking the supremum over \(h\) proves
\(q(\Psi(\zeta v))\le |\zeta|\). Choose a Lempert disk in direction
\(v\) with derivative \(v\) at zero and holomorphic left inverse
\(\ell_v:D\to\mathbb D\). Then \(\ell_v(p)=0\) and
\(d\ell_v(p)v=1\). Because \(D\Psi(0)=\mathrm{Id}\), Schwarz's equality
case gives \(\ell_v(\Psi(\zeta v))=\zeta\). Consequently, we have
\[q(\Psi(\zeta v))=|\zeta|,\] and hence \(q\circ\Psi=\mu\) on \(I\).

This identity makes \(\Psi\) proper. Indeed, if \(C\subset D\) is
compact, then \(\Psi^{-1}(C)\) is closed in \(I\) and is contained in
\(r\overline I\) for some \(r<1\). The only preimage of \(p\) is zero,
where \(D\Psi(0)=\mathrm{Id}\). It follows that the topological degree
of the proper map \(\Psi:I\to D\) is one. Each fiber is a compact analytic
set in a domain of \(\mathbb C^n\), and is thus finite by the maximum
principle. Its local holomorphic multiplicities are positive and sum to
the degree. Therefore every fiber consists of exactly one point of
multiplicity one, at which the derivative is invertible. The holomorphic
inverse function theorem completes the proof.
\end{proof}

We now present the proof of Theorem~\ref{thm1}.

\begin{proof}[Proof of Theorem~\ref{thm1}]
If (i) or (ii) holds, Proposition~\ref{prop1} gives the common extremal.
Lemma~\ref{lem3} gives the local radial field, and Lemma~\ref{lem4}
and Proposition~\ref{prop2} give (iii). The indicatrix is balanced,
so (iii) implies (iv).

Suppose $\Theta:D\to E$ is a biholomorphism onto a balanced domain
with $\Theta(p)=0$. The infinitesimal Kobayashi metric of $E$ at $0$
is bounded below by a positive multiple of the Euclidean norm,
by Lemma~\ref{lem1} and biholomorphic invariance. For $z\in E$,
the disk $\zeta\mapsto\zeta z$ gives $\kappa_E(0;z)\le1$.
Thus $E$ is bounded. It is pseudoconvex by biholomorphic invariance.
For a bounded balanced pseudoconvex domain, the logarithm of its
Minkowski functional is plurisubharmonic and is its Green function
with pole at the center. The lower bound follows by using this
logarithm as a Green competitor. The reverse bound follows by
restriction to each radial disk (see \cite{JP2013,K1991}). Hence
\[I_E^A(0)=E.\] Rotation averaging gives \[K_E(0,0)=\frac{1}{\Vol(E)}.\]
Under a biholomorphism the volume of the indicatrix transforms
by the squared modulus of the complex Jacobian determinant,
and the diagonal kernel transforms by its reciprocal when pulled
in the same direction. Therefore, the product is invariant,
and (iv) implies (i). The equality of the endpoint values in
\eqref{eq:limits} makes $F_p$ constant. This yields (ii).
The asserted strict increase in the other case follows from
monotonicity and the equivalence with (ii).
\end{proof}

\section{Curvature on smooth strongly convex domains}\label{sec5}

In this section, we give the proof of Theorem~\ref{thm2}.

\begin{proof}[Proof of Theorem~\ref{thm2}]
Throughout the proof, $D\subset\mathbb C^n$ is bounded,
$C^\infty$-smooth and strongly convex, and $p\in D$.
We retain the notation of Theorem~\ref{thm2}. 
For convenience, we divide the argument into six steps.

\smallskip
\emph{Step 1: 
}
As $1\in A^2(D)$, we have $K_D(p,p)>0$. By Lemma~\ref{lem1}
and its proof, $I$ is the bounded convex balanced unit ball of the
infinitesimal Kobayashi metric at $p$, and its Minkowski functional
$\mu=\mu_I$ is a complex norm equivalent to the Euclidean norm.
For each $v\in\partial I$, let $\phi_v:\mathbb D\to D$ be the
unique complex geodesic with initial data
\begin{equation}\label{eq:lempert-radial-disk}
 \phi_v(0)=p,\quad \phi_v'(0)=v,\quad
 G(\phi_v(\zeta))=\log|\zeta|.
\end{equation}
The Green identity follows by pulling back Green competitors to
$\mathbb D$ and applying the subharmonic Schwarz lemma, with the
reverse inequality given by the holomorphic left inverse of $\phi_v$.
Uniqueness of the normalized disks gives
\[
 \phi_{e^{i\theta}v}(\zeta)=\phi_v(e^{i\theta}\zeta).
\]
The circular representation $\Phi_p(\zeta v)=\phi_v(\zeta)$ is a
homeomorphism $I\to D$ and a smooth diffeomorphism away from the
centers \cite[Theorem~3]{L1981} (see also \cite{P1983,P1984}).

The Lempert exhaustion $\tau=e^{2G}$, with $\tau(p)=0$, is
continuous, has compact sublevel closures below $1$, and is smooth
and strictly plurisubharmonic off $p$ \cite[Theorem~2.1]{P1984}.
In dimension one these facts also follow from a Riemann map.
By \eqref{eq:lempert-radial-disk}, we have
\begin{equation}\label{eq:lempert-circular-sublevels}
 \tau\circ\Phi_p=\mu^2,\quad \Phi_p(e^sI)=D_s.
\end{equation}
The complex radial field $X$ on $D\setminus\{p\}$ is defined by
$X(\phi_v(\zeta))=\zeta\phi_v'(\zeta)$ for $\zeta\ne0$.
It is well-defined by the phase relation and is smooth by the smooth
diffeomorphism property of $\Phi_p$. Differentiation of
\eqref{eq:lempert-radial-disk} gives $\partial G(X)=1/2$.
Write $\mathcal L_h$ for the Levi form of a smooth real function
$h$, without the factor $i$. The form $\mathcal L_G$ is semipositive,
and $G$ is harmonic on each punctured disk. Thus
$\mathcal L_G(X,X)=0$, and Cauchy--Schwarz for $\mathcal L_G$ gives
$\mathcal L_G(X,Z)=0$ for every $(1,0)$ vector $Z$.
The chain rule gives
\begin{equation}\label{eq:radial-tau-identities}
 \partial\tau(X)=\tau,\quad
 \mathcal L_\tau(X,Z)=\bar\partial\tau(\overline Z).
\end{equation}
In particular, every positive level of $\tau$ is noncritical.

\smallskip
\emph{Step 2: 
}
Fix a small parameter disk $U\subset\mathbb C$ whose real parts
lie in a compact subinterval of $(-\infty,0)$. Set
$s=\operatorname{Re}t$, $c(t)=e^{2s}$, and
\[
 \mathcal E=\{(t,z)\in U\times D:\tau(z)<c(t)\}.
\]
The closure of $\mathcal E$ in $U\times\mathbb C^n$ is proper over
$U$, and its fibers contain $p$, by compactness of the Green
sublevels. The total space is pseudoconvex. An exhaustion is obtained
by adding a plurisubharmonic exhaustion of $U$ to
$-\log(1-e^{G-\operatorname{Re}t})$. This last function is
plurisubharmonic by convex monotone composition and has value zero
at the pole. The sum tends to infinity at the moving boundary and
over the boundary of $U$.

We construct a smooth defining function to apply Wang's formula.
There is $b>0$ with $\tau(z)\ge b|z-p|^2$, by comparison with
the Green function of a Euclidean ball centered at $p$ containing
$D$. On a small ball about $p$, take a smooth regularized maximum
of $\tau$ and $\varepsilon+(b/2)|z-p|^2$, with smoothing width
less than $\varepsilon/2$. One may use
$(x+y+\vartheta(x-y))/2$, where $\vartheta$ is a smooth convex
regularization of $|\cdot|$ with $|\vartheta'|\le1$ agreeing with
$|\cdot|$ outside the smoothing interval. For sufficiently small
$\varepsilon$, the maximum equals the quadratic near $p$ and equals
$\tau$ on a collar of the small ball. Where both arguments are used
they are smooth and strictly plurisubharmonic, so the regularized
maximum is as well, by the Levi chain rule. Choose the ball and
$\varepsilon$ so that all the modified values lie below
$\inf_{t\in U}c(t)$. Outside the small ball retain $\tau$.
Compose the resulting function with a smooth nondecreasing scalar
function that is the identity up to a level above
$\sup_{t\in U}c(t)$ and is constant above a larger level below $1$.
It then extends by the constant outside $D$, since the larger
sublevel is relatively compact. Denote the resulting smooth
function by $\widehat\tau$. Then
\[
 \rho(t,z)=\widehat\tau(z)-c(t)
\]
defines $\mathcal E$ in $U\times\mathbb C^n$, is strictly
plurisubharmonic in $z$ near the fiber closures, and has nonzero
fiber gradient on their boundaries. These are Wang's hypotheses
(A1) and (A2) \cite[Section~2]{W2017}.

Let $V$ be the horizontal lift of $\partial_t$ for
$i\partial\bar\partial(-\log(-\rho))$. Near the moving boundary
$\rho=\tau-c$. There, put $W=\partial_t+X$. Since
$c_t=c_{t\bar t}=c$, equation \eqref{eq:radial-tau-identities} gives
\[
 \partial\rho(W)=\rho,\quad
 \mathcal L_\rho(W,Z)=\bar\partial\rho(\overline Z)
 \quad(Z\text{ vertical}),\quad
 \mathcal L_\rho(W,W)=\rho.
\]
In particular, in the interior part of this collar, we have
\[
 \mathcal L_{-\log(-\rho)}(W,Z)
 =-\frac{\mathcal L_\rho(W,Z)}{\rho}
  +\frac{\partial\rho(W)\bar\partial\rho(\overline Z)}{\rho^2}
 =0\quad(Z\text{ vertical}).
\]
Positivity of the vertical metric gives uniqueness of the horizontal
lift, so $V=\partial_t+X$ throughout this collar. On the boundary,
$W$ is complex tangent and Levi-orthogonal to the vertical tangent
space by the displayed identities. It has Levi norm zero there.
Thus the boundary geodesic curvature $\theta(\rho)$ vanishes
\cite[Definition~2.8 and Lemma~2.9]{W2017}.

Use the unweighted Bergman spaces of holomorphic top forms, with
the constant normalization for which a form
$f\,dz_1\wedge\cdots\wedge dz_n$ has squared norm
$\int|f|^2\dd V$. Write
\[
 \Theta_0=dz_1\wedge\cdots\wedge dz_n,\quad
 k(t)=K_{D_s}(p,p),\quad f_t=K_{D_s}(\cdot,p),
\]
\[
 u_t=f_t\Theta_0,\quad
 \omega_t=i\partial_z\bar\partial_z(-\log(-\rho(t,\cdot))).
\]
Evaluation of the coefficient at $p$ is a holomorphic dual section,
represented by $u_t$, and $\|u_t\|^2=k(t)>0$. It is smooth in
Wang's sense by \cite[Proposition~5.4]{W2017}. Indeed, a fixed
smooth radial density of mass one supported in a common ball about
$p$ represents evaluation by the mean-value property. Define
\[
 B_t=\bar\partial_z(V_{\rm vert}\mathbin{\lrcorner}u_t).
\]
Here $V_{\rm vert}$ is the vertical part in the product coordinates.
With the $d\bar z_j$ factor placed first,
$B_t=(\bar\partial_zV_{\rm vert})\mathbin{\lrcorner}u_t$.
This is a $\bar\partial_z$-closed $(n-1,1)$-form.

We verify its $L^2(\omega_t)$ integrability. Fix $t$ and set
$\delta=-\rho(t,\cdot)$. On the collar let $A=\bar\partial_zX$,
regarded as a map from $T^{0,1}$ to $T^{1,0}$. Differentiating
$\partial\tau(X)=\tau$ in an antiholomorphic direction and using
\eqref{eq:radial-tau-identities} gives $\partial\tau\circ A=0$.
The fiber metric has the expression
\begin{equation}\label{eq:fiber-complete-metric}
 \omega_t=\frac{i\partial_z\bar\partial_z\rho}{\delta}
    +\frac{i\partial_z\rho\wedge\bar\partial_z\rho}{\delta^2}.
\end{equation}
It bounds the Euclidean metric below by $C^{-1}\delta^{-1}$, and
its restriction to $\ker\partial\tau$ is bounded above by
$C\delta^{-1}$. Since $A$ is Euclidean bounded on the collar,
these bounds imply that the $\omega_t$ Hilbert--Schmidt norm of
$A$ is uniformly bounded. In metric orthonormal coordinates,
\[
 |A\mathbin{\lrcorner}u_t|_{\omega_t}^2
 \le C_n|A|_{\omega_t}^2|u_t|_{\omega_t}^2.
\]
The density $|u_t|_{\omega_t}^2dV_{\omega_t}$ is a fixed multiple
of $|f_t|^2\dd V$, because the metric determinant cancels for a top
form. This proves integrability on the collar. The remaining region
is compact in the fiber.

Let $H_t$ be the $L^2(\omega_t)$ harmonic projection in bidegree
$(n-1,1)$. The metric is complete because $\psi=-\log\delta$
tends to infinity at the boundary and satisfies
$|\partial_z\psi|_{\omega_t}\le1$ by
\eqref{eq:fiber-complete-metric}. Indeed, the metric dominates
the rank-one form $i\partial_z\psi\wedge\bar\partial_z\psi$, so
every path tending to the boundary has infinite length. Wang's
harmonic projection for the $\partial$-Laplacian thus agrees with
that for the $\bar\partial$-Laplacian by the K\"ahler identities.
With $\eta_t=\bar\partial_tu_t$, the $(0,1)$ derivative in the
Bergman family, Wang's unweighted formula
\cite[Theorem~5.2, Eq.~(5.18)]{W2017} and $\theta(\rho)=0$ give
\[
 k_{t\bar t}=\|\eta_t\|^2+\|H_tB_t\|_{\omega_t}^2.
\]
The Euclidean volume normalization amounts to a constant positive
Hermitian metric on the trivial line bundle and scales all the form
inner products consistently. Wang's identities for a holomorphic
dual section give $D_tu_t=0$ for the $(1,0)$ part of the Chern
connection \cite[Eqs.~(5.3)--(5.5)]{W2017}. Metric compatibility
then gives $k_{\bar t}=(\eta_t,u_t)$, with the inner product linear
in its first argument. Therefore,
\begin{equation}\label{eq:suita-harmonic-lower}
 \begin{aligned}
 \frac14F_p''(s)=(\log k)_{t\bar t}
 &=\frac{k\|\eta_t\|^2-| (\eta_t,u_t)|^2}{k^2}
   +\frac{\|H_tB_t\|_{\omega_t}^2}{k}
 \ge\frac{\|H_tB_t\|_{\omega_t}^2}{k}\ge0.
 \end{aligned}
\end{equation}
The factor $1/4$ comes from $s=\operatorname{Re}t$, and $2ns$
has zero mixed derivative. The variation theory applies around
every real level and gives $F_p\in C^\infty(( -\infty,0))$.

\smallskip
\emph{Step 3: 
}
Suppose $F_p''(s_0)=0$. Fix $t$ with $\operatorname{Re}t=s_0$
and write $\Omega=D_{s_0}$, $r=e^{s_0}$, $\omega=\omega_t$,
$B=B_t$, $f=f_t$, and $u=u_t$.
By \eqref{eq:suita-harmonic-lower}, $H_tB=0$.

We give the closed-range estimate directly in bidegree $(n-1,0)$,
using the bounded primitive as in Gromov's estimate
\cite[Section~1.4.A]{Gr1991}. Put $\beta=-i\partial\psi$.
Then $\omega=\bar\partial\beta$ and $\beta$ is $\omega$-bounded.
Let $L$ be wedging by $\omega$. Pointwise Lefschetz linear algebra
gives $|Lh|_\omega=|h|_\omega$ in bidegree $(n-1,0)$ in the
standard K\"ahler normalization. For any smooth compactly supported
form $h$ of this type, we have
\[
 Lh=\bar\partial(\beta\wedge h)+\beta\wedge\bar\partial h.
\]
Pairing with $Lh$ and integrating by parts, the K\"ahler commutator
identity gives $\|\bar\partial^*Lh\|_\omega=\|\partial h\|_\omega$,
since $\bar\partial^*h=0$ in bidegree $(n-1,0)$.
The equality of the K\"ahler Laplacians gives
\[
 \|\partial h\|_\omega^2+\|\partial^*h\|_\omega^2
 =\|\bar\partial h\|_\omega^2.
\]
Boundedness of $\beta$ in the two pairings now yields
\[
 \|h\|_\omega^2\le C\|h\|_\omega\|\bar\partial h\|_\omega.
\]
The estimate extends to the maximal $L^2$ domain of $\bar\partial$
in this bidegree. To see this, add a constant $C_0$ so that
$\psi+C_0\ge0$, and use cutoffs equal to one on
$\{\psi+C_0\le R\}$ and zero on $\{\psi+C_0\ge2R\}$, with
$\omega$-derivatives $O(R^{-1})$. Multiplication by these cutoffs
gives convergence in the graph norm. For compactly supported forms,
coordinate mollification gives smooth approximation in the same
graph norm, as $\bar\partial$ acts coefficientwise in the Euclidean
coordinates and the metric is smooth on compact subsets.

Consequently, the maximal operator
\[
 T=\bar\partial:L^2_{n-1,0}(\omega)\longrightarrow
                     L^2_{n-1,1}(\omega)
\]
has closed range. Indeed, the estimate makes a sequence $h_j$
Cauchy when $Th_j$ is Cauchy, and closedness of $T$ provides a
preimage of the limit. The difference between $B$ and its orthogonal
projection onto $\operatorname{ran}T$ lies in
$\ker\bar\partial\cap\ker T^*$, since $B$ is closed and
$\bar\partial T=0$. This difference is its harmonic projection,
which is zero. Hence there is a form $v\in L^2_{n-1,0}(\omega)$
with $\bar\partial v=B$ in distributions.

The metric satisfies
$\omega\ge C^{-1}\delta^{-1}\omega_{\rm Eucl}$ on $\Omega$.
Diagonalize it in a Euclidean unitary frame with eigenvalues
$\lambda_j$. In the density $|v|_\omega^2dV_\omega$, the
coefficient on the wedge omitting the $j$th covector contributes
its squared Euclidean modulus times $\lambda_j\dd V$, up to a
fixed normalization. As $\lambda_j\ge C^{-1}\delta^{-1}$, we have
\begin{equation}\label{eq:suita-zero-trace-weight}
 \int_\Omega\frac{|v|_{\rm Eucl}^2}{\delta}\dd V<\infty.
\end{equation}
The distribution $V_{\rm vert}\mathbin{\lrcorner}u-v$ is
$\bar\partial$-closed, hence a holomorphic $(n-1,0)$-form by the
distributional Cauchy--Riemann equations. Write it as
$Y\mathbin{\lrcorner}\Theta_0$ with $Y$ a holomorphic vector field
on $\Omega$. In the collar where $V_{\rm vert}=X$, we have
\[
 v=(fX-Y)\mathbin{\lrcorner}\Theta_0
\]
almost everywhere. By \eqref{eq:suita-zero-trace-weight}, for
small $\varepsilon>0$, we get
\begin{equation}\label{eq:radial-boundary-integral}
 \int_{\{z\in\Omega:r-\varepsilon<\sqrt{\tau(z)}<r\}}
 \frac{|Y-fX|^2}{r^2-\tau}\dd V<\infty.
\end{equation}
For $v_0\in\partial I$, the vector-valued function
\[
 H_{v_0}(\zeta)=Y(\phi_{v_0}(\zeta))
 -f(\phi_{v_0}(\zeta))\,\zeta\phi_{v_0}'(\zeta),\quad |\zeta|<r
\]
is holomorphic. For $\xi\in S^{2n-1}$, the Euclidean unit sphere,
put $v_\xi=\xi/\mu(\xi)$ and
\[
 M(a)=\int_{S^{2n-1}}|H_{v_\xi}(a)|^2\dd\sigma(\xi)
\]
for $0<a<r$. Pull \eqref{eq:radial-boundary-integral} back by
$z=\Phi_p(w)$. On $r-\varepsilon\le\mu(w)\le r$, the absolute
real Jacobian of $\Phi_p$ is bounded above and below by positive
constants. For $w=av_\xi$, the Euclidean radial coordinate is
$a/\mu(\xi)$. Thus the polar density is
$a^{2n-1}\mu(\xi)^{-2n}\dd a\,d\sigma(\xi)$. Its factor is
bounded above and below there, and
\[
 r^2-\tau(\Phi_p(w))=r^2-a^2
\]
is comparable to $r-a$. It follows that
\begin{equation}\label{eq:radial-circular-integral}
 \int_{r-\varepsilon}^r\frac{M(a)}{r-a}\dd a<\infty.
\end{equation}
Rotation invariance of $d\sigma$ and the phase relation for the
disks give
\[
 H_{v_{e^{i\theta}\xi}}(a)=H_{v_\xi}(ae^{i\theta}),
\]
so
\[
 M(a)=\int_{S^{2n-1}}\frac1{2\pi}\int_0^{2\pi}
       |H_{v_\xi}(ae^{i\theta})|^2\dd\theta\,d\sigma(\xi).
\]
Each circular mean of the squared norm of a holomorphic function
is nondecreasing in the radius, by its power series. Thus $M$ is
nonnegative and nondecreasing. If $M(a_0)>0$ for any $a_0<r$,
the integral in \eqref{eq:radial-circular-integral} diverges at $r$.
Consequently, $M(a)=0$ throughout $0<a<r$. For fixed $a>0$,
the integrand defining $M(a)$ is continuous in $\xi$ by smoothness
of $\Phi_p$ and $X$ away from the pole. It therefore vanishes at
every $\xi$. The identity theorem on each disk gives
$H_{v_0}\equiv0$ for all $v_0\in\partial I$.

As $f(p)=K_\Omega(p,p)>0$, the field $Y/f$ is holomorphic near
$p$ and extends $X$ there. Denote it by $\widetilde X$. The identity
\[
 \widetilde X(\phi_{v_0}(\zeta))=\zeta\phi_{v_0}'(\zeta)
 \quad(|\zeta|\text{ small})
\]
implies $\widetilde X(p)=0$ and $d\widetilde X_p(v_0)=v_0$ for
each $v_0\in\partial I$. Every nonzero vector is a positive
multiple of such a $v_0$, so $d\widetilde X_p=\mathrm{Id}$.

\smallskip
\emph{Step 4: \textup{(ii)}$\Rightarrow$\textup{(iii)}.}

Under the assumption in Step~3, expand
\[
 \widetilde X(p+w)=w+\sum_{k\ge2}X_k(w),\quad
 \phi_{v_0}(\zeta)=p+\sum_{m\ge1}a_m(v_0)\zeta^m,
\]
where $X_k$ are homogeneous holomorphic polynomials of degree $k$
with values in $\mathbb C^n$, and $a_1(v_0)=v_0$.
Comparing the coefficients of $\zeta^m$ in the disk identity gives,
for $m\ge2$,
\[
 (m-1)a_m(v_0)=[\zeta^m]\sum_{k=2}^{m}
 X_k\!\left(\sum_{j=1}^{m-1}a_j(v_0)\zeta^j\right).
\]
Every product contributing to this coefficient has indices summing
to $m$. Inductively, $a_m$ is the restriction of a homogeneous
holomorphic polynomial $P_m$ of degree $m$, with $P_1(w)=w$.

Choose $M_0$ with $|z-p|\le M_0$ on $D$. For $0<R<1$,
Cauchy's estimate on $|\zeta|=R$ gives
$|P_m(v_0)|\le M_0R^{-m}$ uniformly on $\partial I$.
For each $m$, let $R\uparrow1$ to obtain $|P_m(v_0)|\le M_0$.
By homogeneity, $|P_m(w)|\le M_0\mu(w)^m$ on $I$.
The series $p+\sum_{m\ge1}P_m(w)$ converges locally uniformly
on $I$. For $w=av_0$, $a=\mu(w)<1$, it equals
$\phi_{v_0}(a)=\Phi_p(w)$. Thus $\Phi_p$ is holomorphic on $I$.
Its differential at zero is $\mathrm{Id}$ and its real differential
off zero is invertible by the smooth diffeomorphism property.
The holomorphic inverse function theorem, together with bijectivity,
proves that $\Phi_p$ is biholomorphic. This proves
\textup{(ii)}$\Rightarrow$\textup{(iii)}.

\smallskip
\emph{Step 5: \textup{(iii)}$\Rightarrow$\textup{(i)}.}

By rotation averaging and Jensen's inequality, on a bounded balanced
domain the constant function is the minimum-norm holomorphic
function with value one at zero. The Bergman kernel there is
therefore the reciprocal of the volume. If \textup{(iii)} holds,
then $d\Phi_p(0)=\mathrm{Id}$ by the disk derivatives.
The kernel transformation law and
\eqref{eq:lempert-circular-sublevels} give
\[
 K_{D_s}(p,p)=K_{e^sI}(0,0)=\frac{e^{-2ns}}{\Vol(I)},
 \quad K_D(p,p)=\frac1{\Vol(I)}.
\]
This proves \textup{(iii)}$\Rightarrow$\textup{(i)}.

\smallskip
\emph{Step 6: \textup{(i)}$\Rightarrow$\textup{(ii)}.} 

For the two endpoint values, Proposition~\ref{prop1} gives
\eqref{eq:limits}. The left limit follows from the first-order
Green expansion \eqref{eq:green-tangent} and dilation of Bergman
kernels. The right limit follows by increasing exhaustion of $D$
by $D_s$.

Since $F_p''\ge0$ by \eqref{eq:suita-harmonic-lower}, the secant
slope from $a<s$ to $s$ is bounded above by $F_p'(s)$.
The slope tends to zero as $a\to-\infty$ by \eqref{eq:limits}.
Thus $F_p'\ge0$. The nondecreasing derivative tends to zero at
$-\infty$, since a positive limit would contradict finiteness of
$F_p(-\infty)$ after integration. If \textup{(i)} holds, the two
endpoint values agree. The nondecreasing function $F_p$ is then
constant, so \textup{(ii)} holds. Together with Steps~4 and~5,
this establishes the three equivalences.

If the equivalent conditions fail, \eqref{eq:suita-harmonic-lower}
and the negation of \textup{(ii)} give $F_p''(s)>0$ for every $s<0$.
In all cases,
\[
 F_p'(s)=\int_{-\infty}^s F_p''(u)\dd u,
\]
which gives $F_p'(s)>0$ for every $s<0$ when the conditions fail.
Finally, integrating and using Tonelli's theorem gives
\[
 \begin{aligned}
 F_p(0^-)-F_p(-\infty)
 &=\int_{-\infty}^0\int_{-\infty}^s F_p''(u)\dd u\,ds=\int_{-\infty}^0(-u)F_p''(u)\dd u.
 \end{aligned}
\]
Equation~\eqref{eq:limits} identifies the difference with the left
side of \eqref{eq:defect}. This completes the proof.
\end{proof}

\section{Equality loci}\label{sec6}

The set of equality points is invariant under biholomorphic maps.
Theorem~\ref{thm1} gives a description of its geometry in the convex
case, using the orbit theorem of Kaup and Upmeier \cite{KU1976} for the unit ball of
a complex normed space.

\begin{corollary}\label{cor1}
Let $D\subset\mathbb C^n$ be convex and suppose that $p$ is a Suita
equality point with $K_D(p,p)>0$. Let $\Psi:I_D^A(p)\to D$ be as in
Theorem~\ref{thm1}. There is a complex linear subspace
$M\subset\mathbb C^n$ such that the set of all equality points in $D$
is $\Psi(I_D^A(p)\cap M)$. It is also the orbit of $p$ under
$\operatorname{Aut}(D)$. If it has positive dimension, it is
biholomorphic to a bounded symmetric domain in that dimension.
\end{corollary}

\begin{proof}
Set $I=I_D^A(p)$. Every point in the $\operatorname{Aut}(I)$-orbit
of zero is an equality point by invariance. Conversely, if $q\in I$
is an equality point, Theorem~\ref{thm1} gives a biholomorphism from
$I$ to the convex balanced domain $I_I^A(q)$ that sends $q$ to zero.
Both domains are unit balls of complex norms. By the linear
equivalence theorem of Kaup and Upmeier \cite{KU1976}, there is
a complex linear isomorphism taking $I_I^A(q)$ onto $I$.
The composite is an automorphism of $I$ taking $q$ to zero.
The orbit description in \cite{KU1976} gives the subspace $M$.

The orbit is a connected bounded domain in $M$ because $I\cap M$
is convex and relatively open in $M$. The involution \(w\mapsto-w\)
of $I$ restricts to an involution of this orbit with sole fixed
point $0$. Conjugating it by an automorphism of $I$ taking $0$ to
a given point of the orbit yields a holomorphic symmetry at that
point. This is the defining property of a bounded symmetric
domain; see \cite{JP2013}.
\end{proof}

\begin{corollary}\label{cor2}
For a bounded $C^\infty$-smooth strongly convex domain, the set of
Suita equality points is empty, a singleton, or the whole domain.
The last case occurs precisely for domains biholomorphic to the ball.
For every $n\ge2$, the singleton case occurs in $\mathbb C^n$.
\end{corollary}

\begin{proof}
At an equality point $p$, conjugate the rotations
$w\mapsto e^{i\theta}w$ on the indicatrix by $\Psi$.
They give a one-parameter group of automorphisms fixing $p$ and
having differential $e^{i\theta}\mathrm{Id}$ there. Two distinct
points with such groups force $D$ to be biholomorphic to the ball
by Patrizio \cite[Theorem~7]{P1983}. The ball is homogeneous and
balanced at zero, so it has equality at every point.

For the singleton case, let $\epsilon>0$ and put
\[
 E_\epsilon=\{z\in\C^n:|z|^2+\epsilon\textstyle\sum_{j=1}^n|z_j|^4<1\}.
\]
It is smooth, strongly convex, and balanced at zero. If it were
biholomorphic to the ball, composing with a ball automorphism would
give a pointed biholomorphism at zero. Such a map is linear by
Cartan's theorem for bounded circular domains \cite{JP2013}.
A linear ball model would make $\mu_{E_\epsilon}^2$ a positive
Hermitian quadratic form. Independent coordinate rotations and
permutations preserve $E_\epsilon$, so this form would be a scalar
multiple of $|z|^2$. Along a coordinate axis, the boundary radius obeys
\[r^{2}+\epsilon r^{4}=1.\]
For the unit vector \((1,1,0,\dots,0)/\sqrt{2}\), the corresponding boundary radius satisfies
\[r^{2}+\frac{\epsilon}{2}\,r^{4}=1.\]
Since these radii differ, \(E_{\epsilon}\) possesses exactly one equality point.
\end{proof}

\section*{Final remarks}\label{sec7}

In this paper, we resolve the Blocki–Zwonek rigidity conjecture for convex domains in \(\mathbb{C}^n\), including the unbounded case. Our analysis relies on convex-geometric properties, such as the existence of Lempert complex geodesics and compact Green sublevel sets, so the main results cannot be directly extended to general pseudoconvex domains. We construct a holomorphic radial vector field via weighted Bergman kernel estimates, and this new tool is promising for other rigidity problems concerning Bergman kernels and volume-comparison inequalities.

Several open problems remain for further investigation. First, it is interesting to study the Suita equality for finite-type convex domains with non-smooth boundaries. Second, one may ask whether there exist convex domains for which the Suita equality locus has positive intermediate dimension, yet the domain itself is not a bounded symmetric domain. Third, it is natural to explore how much of our analysis can be transplanted to circular or balanced non-convex pseudoconvex domains.

	\vskip.05in
	\noindent{\bf Acknowledgements.} 
  The present investigation was partially supported by the \textit{Key Project of Education Department of Hunan Province} under Grant no. 25A0668, the \textit{Natural Science Foundation of Changsha} under Grant no. kq2502003 of the P. R. China.

	

	\vskip.05in
	\noindent{\bf Conflicts of interest.} The author declares that there is no conflict of interest.
	\vskip.05in
	\noindent{\bf Data availability.} Data sharing does not apply to this article, as no datasets were generated or analyzed during the current study.
	
\vskip .05in
\noindent{\bf AI use disclosure.} Generative AI was used to assist with English language polishing.
All AI-revised text was checked by the author. The author independently carried out all mathematical work and take full responsibility for the manuscript.


\end{document}